\documentclass[12pt,leqno]{article}
\usepackage{amsfonts}
\usepackage{amsmath, amsthm, amsfonts, amssymb, color}
\usepackage{mathrsfs}
\usepackage{color}
\theoremstyle{definition}

\newcommand{\scr}[1]{\mathscr #1}
\definecolor{wco}{rgb}{0.5,0.2,0.3}

\numberwithin{equation}{section} \theoremstyle{remark}

\newcommand{\ua}{\uparrow}

\title{{\bf Gradient and  Hessian Estimates for Dirichlet and Neumann Eigenfunctions}\footnote{Supported in
 part by  NNSFC (11771326, 11431014).} }
\author{
{\bf     Feng-Yu Wang  }\\
\footnotesize{ Center for Applied Mathematics, Tianjin University, Tianjin 300072, China}\\
 \footnotesize{ Department of Mathematics,
Swansea University, Singleton Park, SA2 8PP, United Kingdom}\\
\footnotesize{  wangfy@tju.edu.cn}}
\begin{document}
\allowdisplaybreaks
\def\R{\mathbb R}  \def\ff{\frac} \def\ss{\sqrt} \def\B{\mathbf
B}
\def\N{\mathbb N} \def\kk{\kappa} \def\m{{\bf m}}
\def\ee{\varepsilon}\def\ddd{D^*}
\def\dd{\delta} \def\DD{\Delta} \def\vv{\varepsilon} \def\rr{\rho}
\def\<{\langle} \def\>{\rangle} \def\GG{\Gamma} \def\gg{\gamma}
  \def\nn{\nabla} \def\pp{\partial} \def\E{\mathbb E}
\def\d{\text{\rm{d}}} \def\bb{\beta} \def\aa{\alpha} \def\D{\scr D}
  \def\si{\sigma} \def\ess{\text{\rm{ess}}}
\def\beg{\begin} \def\beq{\begin{equation}}  \def\F{\scr F}
\def\Ric{\mathcal Ric} \def\Hess{\text{\rm{Hess}}}
\def\e{\text{\rm{e}}} \def\ua{\underline a} \def\OO{\Omega}  \def\oo{\omega}
 \def\tt{\tilde}
\def\cut{\text{\rm{cut}}} \def\P{\mathbb P} \def\ifn{I_n(f^{\bigotimes n})}
\def\C{\scr C}      \def\aaa{\mathbf{r}}     \def\r{r}
\def\gap{\text{\rm{gap}}} \def\prr{\pi_{{\bf m},\varrho}}  \def\r{\mathbf r}
\def\Z{\mathbb Z} \def\vrr{\varrho} \def\ll{\lambda}
\def\L{\scr L}\def\Tt{\tt} \def\TT{\tt}\def\II{\mathbb I}
\def\i{{\rm in}}\def\Sect{{\rm Sect}}  \def\H{\mathbb H}
\def\M{\scr M}\def\Q{\mathbb Q} \def\texto{\text{o}} \def\LL{\Lambda}
\def\Rank{{\rm Rank}} \def\B{\scr B} \def\i{{\rm i}} \def\HR{\hat{\R}^d}
\def\to{\rightarrow}\def\l{\ell}\def\iint{\int}
\def\EE{\scr E}\def\Cut{{\rm Cut}}
\def\A{\scr A} \def\Lip{{\rm Lip}}\def\S{\mathbb S}
\def\BB{\scr B}\def\Ent{{\rm Ent}} \def\i{{\rm i}}\def\itparallel{{\it\parallel}}
\def\g{{\mathbf g}}\def\Sect{{\mathcal Sec}}\def\T{\mathcal T}\def\Eig{{\rm Eig}}\def\Sect{{\rm Sect}}
\maketitle

\begin{abstract} We establish  integral formulas and  sharp two-sided bounds for the Ricci curvature,   mean curvature and   second fundamental form on a Riemannian manifold with boundary. As applications,  sharp gradient and Hessian   estimates   are derived  for the  Dirichlet and Neumann eigenfunctions.
\end{abstract} \noindent
 \vskip 2cm

\section{Introduction}

Let $(M,\g) $ be a $d$-dimensional complete connected Riemannian manifold with boundary $\pp M$, and let $N$ be the inward unit normal vector field of $\pp M$.   We also denote $\g(u,v)= \<u,v\>$ for two vector fields $u,v$.
For $V\in C^2(M)$, let  $L = \DD+\nn V$ and $\mu(\d x)=\e^{V(x)}{\rm vol}(\d x),$ where {\rm vol} is the volume measure. Then $L$ is symmetric in $L^2(\mu)$ under the Neumann condition ($Nf|_{\pp M}=0$) or the Dirichlet condition ($f|_{\pp M}=0$). We estimate the gradient and Hessian of the Dirichlet and Neumann eigenfunctions for $L$ by using the following quantities:
\beg{enumerate} \item[$\bullet$] Bakry-Emery curvature on $M$: $\Ric^V= \Ric-\Hess_V$, where $\Ric$ is the Ricci curvature.
\item[$\bullet$] Second fundamental form of $\pp M$: $\II_{\pp}(u,v)=-\<\nn_u N,v\>,\ \ u,v\in T\pp M$.
\item[$\bullet$] Weighted mean curvature of $\pp M$: $H_V ={\rm tr} (\II_{\pp})- NV = -L \rr_\pp$ on $\pp M$, where $\rr_\pp$ is the Riemannian distance to the boundary $\pp M$. When $V=0$, $H_0$ is the usual mean curvature of $\pp M$. \end{enumerate}

Let $\mu_\pp$ be the area measure on $\pp M$ induced by $\mu$.   For a nonnegative function $f$, we denote
$$\mu(f)=\int_Mf\d\mu,\ \ \mu_\pp(f)= \int_{\pp M}f\d\mu_\pp.$$
We call $(\ll,\phi)\in (0,\infty)\times C^2(M)$ an eigenpair of $L$, if $L\phi=-\ll\phi$ holds.
Let $\Eig_N(L)$ be the set of eigenpairs $(\ll, \phi)$ with $\mu(\phi^2)=1$ for the Neumann problem (i.e. $N\phi|_{\pp M}=0$), and let $\Eig_D(L)$ be that for the Dirichlet eigenproblem (i.e. $\phi|_{\pp M}=0$). We aim to estimate
$\mu_\pp(\phi^2),\ \mu_\pp(|\nn \phi|^2)$ and $\mu (\|\Hess_\phi\|_{HS}^2)$ for Dirichlet and Neumann eigenfunctions respectively. Before state our results, we first recall known boundary estimates derived in \cite{BHT,HT}.

According to \cite{HT}, there exists a constant $C>0$ such that
\beq\label{4C1} \mu_\pp(|\nn \phi|^2) \le C\ll,\ \ (\ll,\phi)\in \Eig_D(\DD).\end{equation}
When $M$ does not contain any trapped geodesic, i.e. any geodesic starting from a point in $M$ will eventually go beyond $M$ (it is the case when $M$ is a domain in $\R^d$), then there exists a constant $c>0$ such that
\beq\label{4C1'}\mu_\pp(|\nn \phi|^2) \ge c\ll,\ \ (\ll,\phi)\in \Eig_N(\DD).\end{equation} So, as a general result, the order of $\ll$ in \eqref{4C1} is sharp. But, in general, \eqref{4C1'} is not true, see \cite{HT} for counterexamples, which include semi-spheres and cylinders. This indicate that for   boundary estimates of eigenfunctions, a smooth domain  in $\R^d$ may be essentially different from a Riemannian manifold with boundary.

However,   for Neumann eigenfunctions the estimate \eqref{4C1} does not hold.  According to \cite{BHT}, when $M$ is a bounded smooth domain in $\R^d$ and $L=\DD$, there exists a constant $C>0$ such that
\beq\label{4C2} \mu_\pp(\phi^2)\le C \ll^{\ff 1 3},\ \ \mu_\pp(|\nn \phi|^2) \le C\ll^{\ff 4 3},\ \ (\ll,\phi)\in \Eig_N(\DD),\end{equation}
where the order of $\ll$ in both estimates is sharp for the disc in $\R^2$.

In this paper, we aim to derive sharp Hessian estimate for $(\ll,\phi)\in \Eig_D(L)$, and extend \eqref{4C2} to general compact Riemannian manifolds with boundary  which, in turn, to imply sharp Hessian estimates for $(\ll,\phi)\in\Eig_N(L)$.

We first consider the Hessian estimate for  Dirichlet eigenfunctions.  We will see that the following result is a  straightforward consequence of
\eqref{4C1} and the integral formula \eqref{IF2} proved in the next section.

\beg{thm}[Diriclet eigenfunctions]\label{T4.1}  Let $M$ be a $d$-dimensional connected compact Riemannian manifold with boundary $\pp M$.   Let $K_1,K_2,\dd_1,\dd_2\in \R$ be constants such that
$$K_1\le\Ric^V\le K_2,\ \ \dd_1\le H_V\le \dd_2.$$
Then
\beq\label{4.1} K_1\ll +\dd_1\mu_\pp(|\nn \phi|^2) \le \ll^2-\mu (\|\Hess_\phi\|_{HS}^2)\le K_2\ll+ \dd_2 \mu_\pp(|\nn\phi|^2),\ \ (\ll,\phi)\in \Eig_D(L).\end{equation}
Consequently, there exists a constant $C>0$ such that
\beq\label{4.10} |\mu (\|\Hess_\phi\|_{HS}^2)-\ll^2|\le C\ll,\ \ (\ll,\phi)\in\Eig_D(L).\end{equation} In particular, if $\Ric^V=0$ and $H_V=0$, then $\eqref{4.10}$ holds for $C=0$.  \end{thm}
\beg{proof} Since $L\phi=-\ll \phi$, \eqref{4.1}   follows from \eqref{IF2} in the next section. Next, by repeating the argument in \cite{HT}, we may prove \eqref{4C1} for $\Eig_D(L)$ replacing $\Eig_D(\DD)$. So, \eqref{4.10} follows from \eqref{4.1}.
\end{proof}
By \eqref{4.1} and the sharpness of \eqref{4C1} as explained above, the order of $\ll$ in estimate \eqref{4.10} is sharp as well.

The situation for  the Neumann problem is more complicated. We address the main result below but leave the proof to Section 3.

\beg{thm}[Neumann eigenfunctions]\label{T4.2} Let $M$ be a $d$-dimensional connected compact Riemannian manifold with boundary $\pp M$. Let $K_1,K_2,\kk_1,\kk_2\in \R$ be constants such that
$$K_1\le\Ric^V\le K_2,\ \ \kk_1\le \II_\pp\le \kk_2.$$ Then \beq\label{4.3} K_1\ll +\kk_1\mu_\pp(|\nn \phi|^2) \le \ll^2-\mu(\|\Hess_\phi\|_{HS}^2)\le K_2\ll+ \kk_2 \mu_\pp(|\nn\phi|^2),\ \ (\ll,\phi)\in \Eig_N(L).\end{equation}
Moreover, there exists a constant $C>0$ such that
\beq\label{4C20} \mu_\pp( \phi^2) \le C\ll^{\ff 1 3},\ \ \mu_\pp(|\nn \phi|^2) \le C\ll^{\ff 4 3},\ \ (\ll,\phi)\in \Eig_N(L).\end{equation}
Consequently, there exists a constant $C>0$ such that
\beq\label{4.2} |\mu (\|\Hess_\phi\|_{HS}^2)-\ll^2|\le C\ll^{\ff 4 3},\ \ (\ll,\phi)\in\Eig_N(L).\end{equation} If in particular $\Ric^V=0$ and $\II_\pp=0$, $\eqref{4.2}$ holds for $C=0.$
\end{thm}

By \eqref{4.3} and the sharpness of \eqref{4C2} for   the disc in $\R^2$ as explained in \cite{BHT}, the order of $\ll$ in \eqref{4.2} is sharp as well.

In Section 2, we establish integral formulas and two-sided bounds for the above mentioned geometry quantities, which   will be   used   in Section 3 to prove   Theorem \ref{T4.2}.

 \section{Integral characterizations  of $\Ric^V, \II_{\pp}$ and $H_V$}

In this section, we assume that $M$ is a  Riemannian manifold with boundary which is not necessarily compact nor connected.
Let $C_0^\infty(M)$ be the set of smooth functions on $M$ with compact support. We
consider the  following two classes of reference functions for the Neumann and Dirichlet problems respectively:
\beg{align*} C_N^\infty(M):=\big\{f\in C_0^\infty(M):\ Nf|_{\pp M}=0\big\},\ \
 C_D^\infty(M):=\big\{f\in C_0^\infty(M):\ f|_{\pp M}=0\big\}.\end{align*}
  By Bochner-Weizenb\"ock  and integration by parts formulas, we have the following integral formulas for $\Ric^V, \II_\pp$ and $H_V$.

\beg{thm}\label{T1} Let $C_N^\infty(M)$ and $C_D^\infty(M)$ be in above. We have
\beq\label{IF1} \int_M\Big\{(Lf)^2-\|\Hess_f\|_{HS}^2 -\Ric^V(\nn f,\nn f)\Big\}\,\d\mu= \int_{\pp M} \II_\pp(\nn f,\nn f)\,\d\mu_\pp,\ \ f\in C_N^\infty(M);\end{equation}
\beq\label{IF2}  \int_M\Big\{(Lf)^2-\|\Hess_f\|_{HS}^2 -\Ric^V(\nn f,\nn f)\Big\}\,\d\mu= \int_{\pp M} H_V |\nn f|^2\,\d\mu_\pp,\ \ f\in C_D^\infty(M).\end{equation}\end{thm}

\beg{proof} By Bochner-Weizenb\"ock formula,
\beq\label{BW}\ff 1 2 L |\nn f|^2 = \|\Hess_f\|_{HS}^2 +\<\nn Lf,\nn f\> + \Ric^V(\nn f,\nn f),\ \ f\in C^\infty(M).\end{equation}
Next, the integration by parts formula gives
\beq\label{IT1}   -\int_M\<\nn Lf,\nn f\>\d\mu =\int_{\pp M} (Lf)Nf\,\d\mu_\pp +\int_M (Lf)^2\d \mu, \end{equation}
\beq\label{IT2} \ff 1 2 \int_M L|\nn f|^2\d\mu= - \ff 1 2 \int_{\pp M} N|\nn f|^2\d\mu_\pp=- \int_{\pp M} \Hess_f(\nn f, N) \d\mu_\pp. \end{equation}
Integrating \eqref{BW} with respect to $\d\mu$ and using \eqref{IT1}, \eqref{IT2}, we arrive at
\beq\label{IT3} \beg{split} &\int_M\Big\{(Lf)^2-\|\Hess_f\|_{HS}^2 -\Ric^V(\nn f,\nn f)\Big\}\,\d\mu\\
&=\int_{\pp M} \big\{\Hess_f(N,\nn f)-(Nf)Lf\big\}\,\d\mu_\pp,\ \ f\in C_0^\infty(M).\end{split}\end{equation}
With this formula we are able to prove  \eqref{IF1} and \eqref{IF2} as follows.

Firstly, for $f\in C_N^\infty(M)$, we have $Nf|_{\pp M}=0$ and, by \cite[the formula after (3.2)]{W09},
$$ \Hess_f(N,\nn f)|_{\pp M}=  -\<\nn_{\nn f}N,\nn f\>|_{\pp M}=\II_\pp(\nn f,\nn f)|_{\pp M}.$$
Then    \eqref{IF1} follows from \eqref{IT3}.

Next, for $f\in C_D^\infty(M)$, we have $f|_{\pp M}=0$. So,    $\nn f|_{\pp M}=(Nf)N$ and
\beq\label{HS2} \Hess_f(N,\nn f)|_{\pp M}= (Nf)\Hess_f(N,N)|_{\pp M}.\end{equation}
Let $\{v_i\}_{i=1}^{d-1}$ be orthonormal vector fields in a neighborhood of a point $x\in \pp M$, such that   $\nn v_i(x)=0$ and $\<N,v_i\>(x)=0$. Then
$$ \DD f(x)= \Hess_f(N,N)(x) +\sum_{i=1}^{d-1} \Hess_f(v_i,v_i)(x).$$
Combining this with $\nn f|_{\pp M}= (Nf)N|_{\pp M}$ and $\<v_i,N\>(x)=0=\<N,\nn_{v_i}v_i\>(x)$, we arrive at
\beg{align*} \DD f(x) -\Hess_f(N,N)(x) &=\sum_{i=1}^{d-1} \Hess_f(v_i,v_i)(x)= \sum_{i=1}^{d-1} v_i\<\nn f, v_i\>(x)\\
&= \sum_{i=1}^{d-1} v_i\{(Nf) \<N, v_i\>\}(x)= \sum_{i=1}^{d-1} \{(Nf)\<\nn_{v_i}N,v_i\>\}(x)\\
&=-\{(Nf){\rm tr}(\II_\pp)\}(x)= -( H_0 Nf)(x).\end{align*} This, together with \eqref{HS2}, yields
$$\Hess_f(N,\nn f)|_{\pp M} = (Nf) (\DD f+ H_0 Nf)= \big\{ H_0 (Nf)^2+ (Nf)\DD f\big\}|_{\pp M}.$$ Combining with     $\nn f|_{\pp M}= (Nf)N|_{\pp M}$  leads to
\beg{align*} & \{\Hess_f(N,\nn f) -(Nf) Lf\}|_{\pp M} =\{H_0 |\nn f|^2 -(Nf)\<\nn V,\nn f\>\}|_{\pp M}\\
&= |\nn f|^2 (H_0-NV)|_{\pp M}= H_V |\nn f|^2|_{\pp M}.\end{align*} Substituting into   \eqref{IT3}, we prove   \eqref{IF2}. \end{proof}

We now   characterize bounds of $\Ric^V,\ \II_\pp$ and $H_V$. For a symmetric $2$-tensor $Q$, we write $Q\ge 0$ (or $Q\le 0$)  if $Q(v,v)\le 0$ (or $Q(v,v)\le 0$) holds for all vectors $v$. For two symmetric $2$-tensors $Q_1,Q_2$, we write $Q_1\ge Q_2$ (equivalently, $Q_2\le Q_1$) if $Q_1-Q_2\ge 0$ (equivalently, $Q_2-Q_1\le 0$).

 \beg{thm}\label{T2} Let $Q$ and $Q_\pp$ be continuous  symmetric 2-tensors on $M$ and $\pp M$ respectively, and let $q\in C(\pp M)$.
 \beg{enumerate} \item[$(1)$] $\Ric^V\ge Q$ and $\II_{\pp}\ge Q_\pp$ if and only if
 \beq\label{E1} \beg{split}&\int_{M} \big\{(L f)^2- \|\Hess_f\|_{HS}^2- Q(\nn f,\nn f)\big\}\d\mu \\
 &\ge \int_{\pp M} Q_\pp (\nn f,\nn f) \d\mu_\pp,\ \ f\in C_{N}^\infty(M).\end{split}\end{equation}
  \item[$(2)$] $\Ric^V\le Q$ and $\II_{\pp}\le Q_\pp$ if and only if
 \beq\label{E2}\beg{split} & \int_{M} \big\{(L f)^2- \|\Hess_f\|_{HS}^2- Q(\nn f,\nn f)\big\}\d\mu \\
 &\le \int_{\pp M} Q_\pp (\nn f,\nn f) \d\mu_\pp,\ \ f\in C_{N}^\infty(M).\end{split}\end{equation}
\item[$(3)$] $\Ric^V\ge Q$ and $H_V\ge q$ if and only if
 \beq\label{E1'} \beg{split}&\int_{M} \big\{(L f)^2- \|\Hess_f\|_{HS}^2- Q(\nn f,\nn f)\big\}\d\mu \\
 &\ge \int_{\pp M} q|\nn f|^2 \d\mu_\pp,\ \ f\in C_{D}^\infty(M).\end{split}\end{equation}
  \item[$(4)$] $\Ric^V\le Q$ and $H_V\le q$ if and only if
 \beq\label{E2'}\beg{split} & \int_{M} \big\{(L f)^2- \|\Hess_f\|_{HS}^2- Q(\nn f,\nn f)\big\}\d\mu \\
 &\le \int_{\pp M} q|\nn f|^2  \d\mu_\pp,\ \ f\in C_{D}^\infty(M).\end{split}\end{equation}
\end{enumerate} \end{thm}

\beg{proof} According to Theorem \ref{T1}, we only need to prove the sufficiency in all assertions.

According to \eqref{IF1} and \eqref{IF2},  the inequalities  \eqref{E1}  and \eqref{E1'}  are equivalent to the following ones respectively:
\beg{align*} & \int_M\big\{\Ric^V-Q\big\}(\nn f,\nn f)\d\mu + \int_{\pp M} \big\{\II_\pp- Q_\pp\big\}(\nn f,\nn f)\d\mu_\pp\ge 0,\ \ f\in C_N^\infty(M),\\
&\int_M\big\{\Ric^V-Q\big\}(\nn f,\nn f)\d\mu + \int_{\pp M} \big\{(H_V-q)|\nn f|^2\big\}\d\mu_\pp\ge 0,\ \ f\in C_D^\infty(M).\end{align*}
By the following Lemma  \ref{LP2}, the first implies $\Ric^V\ge Q$ and $\II_\pp\ge Q_\pp$, while the second yields $\Ric^V\ge Q$ and   $H_V\ge q$. Thus, assertions (1) and (3) hold.  Similarly, we can prove assertions (2) and (4). \end{proof}

\beg{lem}\label{LP2} Let $Q,Q_\pp$ be continuous symmetric $2$-tensors on $TM$ and $T\pp M$ respectively, and let $h\in C(\pp M)$.
\beg{enumerate} \item[$(1)$]   $Q\ge 0$ and $Q_\pp \ge 0$ if and only if
\beq\label{E3} \int_M  Q (\nn f,\nn f) \d\mu +\int_{\pp M}  Q_\pp  (\nn f,\nn f) \d\mu_\pp\ge 0, \ \ f\in C_{N}^\infty(M).\end{equation}
\item[$(2)$] $Q\ge 0$ and $h\ge 0$ if and only if
\beq\label{E4} \int_M  Q (\nn f,\nn f) \d\mu +\int_{\pp M}  h|\nn f|^2 \d\mu_\pp\ge 0, \ \ f\in C_D^\infty(M).\end{equation}\end{enumerate}
  \end{lem}

\beg{proof} The necessity in these assertions are trivial. Below we prove the   sufficiency.

(a) $Q\ge 0$. For $f\in C_0^\infty(M\setminus \pp M)
\subset C_N^\infty(M)\cap C_D^\infty(M)$, we have $\nn f|_{\pp M}=0$ so that each of   \eqref{E3} and \eqref{E4} implies
$$\int_{M} Q(\nn f,\nn f)\d\mu  \ge 0,\ \ f\in C_0^\infty(M\setminus \pp M).$$
According to \cite[Lemma 2.2]{W17} for $M\setminus \pp M$ replacing $M$, this implies $Q\ge 0$   in $M\setminus \pp M$. By the continuity of $Q$,  it holds on $M$.

(b)  $Q_\pp\ge 0$.
 Let $x_0\in \pp M$ and $X_0\in T_{x_0}\pp M$ with $|X_0|=1$, we aim to prove $Q_\pp (X_0,X_0)\ge 0.$ To this end, we take the normal coordinates   in a neighborhood $O(x_0)$ of $x_0$ such that
\beg{enumerate} \item[$(1)$] $x_0=0\in\R^d, \ X_0=\pp_1|_{x=0};$
\item[$(2)$] For some constant $r_0>0$, $$O(x_0)= \bigg\{(x^1,\cdots,x^d)\in \R^d: 0\le x^d,\  \sum_{i=1}^{d} |x^i|^2<r_0\bigg\};$$
\item[$(3)$]   $(\pp M)\cap O(x_0)= \big\{x=(x^1,\cdots,x^d)\in O(x_0): x^d=0\big\},$ on which  $N=\pp_d$. \end{enumerate}
Under this local coordinate system, let $ \hat x = (x^1,\cdots, x^{d-1},0)\ \text{for}\ x=(x^1,\cdots, x^d).$  Then there exist   symmetric matrix-valued continuous functions $(q^{ij})_{1\le i,j\le d}$ and $(q^{ij}_\pp)_{1\le i,j\le d-1}$ such that
\beq\label{E4'} \beg{split} &Q(\nn f,\nn f)\d\mu= \sum_{i,j=1}^d \big\{q^{ij}  (\pp_if)(\pp_jf) \big\}(x)\d x\  \text{on} \ O(x_0),\\
&Q_\pp (\nn f,\nn f)\d\mu_\pp = \sum_{i,j=1}^{d-1} \big\{q^{ij}_\pp  (\pp_if)(\pp_jf) \big\}(\hat x)\d \hat x\    \text{on} \ (\pp M)\cap O(x_0).\end{split}\end{equation}

Now, for any $n\ge 1$ and $x\in \R^d$, let
$$\big(\phi_n(x)\big)^i= \beg{cases} n^2 x^1,\ \ & i= 1,\\
n x^i,\ \ & 2\le i\le d.\end{cases}$$ Let $f\in C_0^\infty(O(x_0))$ with $Nf|_{\pp M}=0$, i.e. $\pp_d f|_{x^d=0}=0$.   Then
$$f_n:= f\circ \phi_n\in C_0^\infty(O(x_0)),\ \ \pp_d f_n|_{x^d=0}=0,\ \ n\ge 1.$$
So, by \eqref{E3} and \eqref{E4'} we obtain
\beg{align*} &0  \le  \int_{\R^d}  \sum_{i,j=1}^d \big\{q^{ij}(\pp_i f_n)(\pp_j f_n)\big\}(x)\d x +\int_{\R^{d-1}} \sum_{i,j=1}^{d-1} \big\{q^{ij}_\pp(\pp_i f_n)(\pp_j f_n)\big\}(\hat x)\d \hat x\\
&=    \int_{\R^d} \bigg\{ n^{1-d} \sum_{i,j=2}^d  (q^{ij}\circ\phi_n^{-1}) (\pp_i f)(\pp_j f)    + 2 n^{2-d} \sum_{j=2}^d  (q^{1j}\circ\phi_n^{-1}) (\pp_1 f)(\pp_j f)\\
&\qquad\qquad  +n^{3-d} (q^{11}\circ\phi_n^{-1}) (\pp_1 f)^2\bigg\}(x)\d x\\
&\quad +  \int_{\R^{d-1}} \bigg\{n^{2-d} \sum_{i,j=2}^{d-1}  (q^{ij}_\pp\circ\phi_n^{-1}) (\pp_i f)(\pp_j f)    + 2 n^{3-d} \sum_{j=2}^{d-1}  (q^{1j}_\pp\circ\phi_n^{-1}) (\pp_1 f)(\pp_j f)\\
 &\qquad\qquad +n^{4-d} (q^{11}_\pp\circ\phi_n^{-1}) (\pp_1 f)^2\bigg\}(\hat x)\d \hat x.\end{align*}
Multiplying by $n^{d-4}$ and letting $n\to \infty$, we arrive  at
$$0\le q^{11}_\pp(0) \int_{\R^{d-1}} (\pp_1 f)^2(\hat x)\d\hat x,\ \ f\in (\pp M)\cap C_0^\infty(O(x_0)).$$ Combining this with the second equality in \eqref{E4'} and noting that $X_0=\pp_1|_{x_0}$, we obtain $Q_\pp(X_0,X_0)\ge 0.$

(c) $h\ge 0$.
Let $g\in C^\infty_0(\pp M)$ with compact support $D\subset \pp M$.    There exist a neighborhood $\scr O$ in $M$ of  $D$, and a constant  $r_0>0$, such that $\rr_\pp\in C^\infty_b(\scr O)$ and the Fermi coordinates
 $$ \scr O\ni x =(\theta, r)\in    \pp M \times [0, r_0)$$ exists, where $x=(\theta, r)$ means $x=\exp_\theta [r N].$
 Let $\gg\in C_0^\infty([0,\infty))$ such that $\gg|_{[0,r_0/2]}=1, \gg|_{[r_0,\infty)}=0.$ For any $n\ge 1$, define
 $$f_n(x):= \beg{cases}  g(\theta) r \gg(n r),  \ &\text{if}\  x=(\theta,r)\in \scr O,\\
 0,\ &\text{otherwise}.\end{cases}$$ Then $f_n\in C_D^\infty(M)$ and
 $$|\nn f_n|^2|_{\pp M}= g^2|_{\pp M},\ \ |\nn f_n|\le c\big\{\|\nn^{\pp M}  g\|_\infty + \|g\|_\infty (1+r_0\|\gg'\|_\infty)\big\}1_{\{\rr_\pp\le r_0/n\}},$$where $\nn^{\pp M}$ is the gradient on $\pp M$.
 So, applying \eqref{E4}  for $f_n$ replacing $f$£¬ we may find out a constant $C>0$ such that for any $n\ge 1$,
 \beg{align*}& \int_{\pp M} (h g^2 )\d\mu_\pp\ge -\int_M Q(\nn f_n,\nn f_n)\d\mu \\
 &\ge - C  \int_{\{\rr_\pp\le   r_0/n\}\cap\scr O} \Big\{ \|\nn^{\pp M}g\|_\infty^2  + \|g\|_\infty^2 (1+r_0\|\gg'\|_\infty)^2\Big\}\d\mu.\end{align*}
 By letting $n\to\infty$ we arrive at
 $$\int_{\pp M} (h g^2 )\d\mu_\pp\ge 0,\ \ g\in C_0^\infty(\pp M),$$
 which implies $h\ge 0$ as $g\in C_0^\infty(\pp M)$ is arbitrary. \end{proof}

\section{Proof of Theorem \ref{T4.2}}

To prove Theorem \ref{T4.2}, we present some lemmas.

\beg{lem}\label{L3} There exists a constant $C>0$ such that
$$\int_{\pp M} |\nn \phi|^2\d\mu_\pp\le C  \bigg(\ll + \ll \int_{\pp M}\phi^2\d\mu_\pp\bigg),\ \ (\ll,\phi)\in \Eig_N(L).$$\end{lem}

\beg{proof}  Let $r_0>0$ such that $\rr_\pp$ is smooth on
$\pp_{r_0}M:= \{\rr_\pp\le r_0\}$  and the Fermi coordinate system $x=(\theta, r)\in  \pp M\times [0,r_0]$ exists on $M_{0,r_0}:=\{\rr_\pp\le r_0\}$. Under this coordinate system we have
\beq\label{LP} L= \aa \{\DD_{\pp M} +\nn^{\pp M}V\}+ \pp_r^2 + Z,\end{equation} where $\aa\in C^\infty(\pp M\times [0,r_0])$ is strictly positive with $\aa(\cdot,0)=1$, $\DD_{\pp M}$ and $\nn^{\pp M}$ are  the Laplacian and gradient on the $(d-1)$-dimensional Riemannian manifold  $\pp M$ respectively, and $Z$ is a $C^1$ (hence, bounded) vector field on $M$. Using the integration by parts formula on $\pp M$,   \eqref{LP}, and $L\phi=-\ll\phi$, we obtain
\beg{align*}&\int_{\pp M} |\nn \phi|^2\d\mu_\pp = -\int_{\pp M}\phi \{ \DD_{\pp}+\nn^{\pp M} V\} \phi \d\mu_\pp \\
&= \int_{\pp M} \Big\{  \phi\Hess_\phi(N,N) +\phi Z\phi- \phi L\phi\Big\}\d\mu_\pp\\
&\le (\ll+\|Z\|_\infty^2)\int_{\pp M}\phi^2\d\mu_\pp +\ff 1 4 \int_{\pp M} |\nn \phi|^2\d\mu_\pp +\int_{\pp M}\phi\Hess_\phi(N,N)\d\mu_\pp.\end{align*} Since  $\ll\ge \ll_1^N>0$,  this implies
\beq\label{LP2} \int_{\pp M}|\nn \phi|^2\d\mu_\pp\le c_1 \ll \int_{\pp M} \phi^2\d\mu_\pp + \ff 4 3 \int_{\pp M}\phi\Hess_\phi(N,N)\d\mu_\pp,\ \ (\ll,\phi)\in\Eig_N(L)\end{equation} for some constant $c_1>0$.
To estimate $\int_{\pp M}\phi\Hess_\phi(N,N)\d\mu_\pp$, we take $\gg\in C_0^\infty([0,\infty))$ such that $\gg|_{[0,r_0/2]}=1, \gg|_{[2r_0/3,\infty)}=0$. By $L\phi=-\ll\phi$, $N\phi|_{\pp M}=0$ and using  integration by parts,  we have
\beq\label{LP3} \beg{split} &\int_{\pp M} \phi \Hess_\phi(N,N)\d\mu_\pp = \int_{\pp M} \phi N\<\gg(\rr_\pp) \nn \rr_\pp, \nn\phi\> \d\mu_\pp\\
&= \int_M\Big\{- \phi L \<\gg(\rr_\pp) \nn \rr_\pp, \nn\phi\> + (L \phi) \<\gg(\rr_\pp) \nn \rr_\pp, \nn\phi\>\Big\}\d\mu\\
&=\int_M  \phi [\gg(\rr_\pp) \nn \rr_\pp,L]\, \phi\d\mu,\end{split} \end{equation}
where $[\gg(\rr_\pp) \nn \rr_\pp,L]:= (\gg(\rr_\pp) \nn \rr_\pp)L - L(\gg(\rr_\pp) \nn \rr_\pp)$ is a continuous  second order differential operator on the compact set $\{\rr_\pp\le r_0\}.$  Combining this with   $\int_M\phi^2\d\mu=1$, we derive
\beq\label{LPQ}\bigg|\int_{\pp M} \phi \Hess_\phi(N,N)\d\mu_\pp \bigg|\le c_2 \bigg(\int_M (\phi^2+|\nn \phi|^2+ \|\Hess_\phi\|_{HS}^2)\d\mu\bigg)^{\ff 1 2}\end{equation}
for some constant $c_2>0$. Combining with \eqref{LP2}, $\mu(\phi^2)=1$ and $\mu(|\nn \phi|^2)=\ll$,  we arrive at
\beq\label{LPP} \int_{\pp M}|\nn \phi|^2\d\mu_\pp\le  c_1 \ll \int_{\pp M} \phi^2\d\mu_\pp + \ff {4c_2} 3 \bigg(1+\ll+\int_{\pp M} \|\Hess_\phi\|_{HS}^2\d\mu\bigg)^{\ff 1 2}.\end{equation}  But by \eqref{4.3} we have
$$\int_{\pp M} \|\Hess_\phi\|_{HS}^2\d\mu\le \ll^2 +c_3 \ll + \int_{\pp M}|\nn \phi|^2\d\mu_\pp$$
for some constant $c_3>0$,     \eqref{LPP} implies the desired estimate for some constant $C>0$.
\end{proof}

\beg{lem}\label{TP} There exists a constant $C>0$ such that
$$\int_{\pp M}\phi^2\d\mu_\pp\le C\ll^{\ff 1 3},\ \ (\ll,\phi)\in \Eig_N(L).$$\end{lem}

We first prove a priori estimate then make improvement. To this end, we introduce some notation.

For any $r\ge 0$, let $\pp_r M= \{\rr_\pp=r\}$ and $\mu_\pp^r$ be the area measure on it induced by $\mu$. For $0<r_1<r_2$, let
$$ M_{r_1,r_2}=\{r_1\le \rr_\pp\le r_2\},\ \ \mu_\pp^{r_1,r_2}= 1_{\pp_{r_1}M} \mu_\pp^{r_1}+  1_{\pp_{r_2}M} \mu_\pp^{r_2}.$$ Obviously,
$\pp_0 M= \pp M, \mu_\pp^0=\mu_\pp.$

Let $\dd>0$ such that $\rr_\pp\in C_b^\infty(M_{0,\dd})$, and   the Fermi coordinate system $(\theta, r)\in \pp M\times [0,\dd]$ gives a diffeomorphism
between $M_{0,\dd}$ and $\pp M\times [0,\dd]$. Under this coordinate system we have
\beq\label{**0}
\mu(\d\theta,\d r)= \psi(\theta, r)\mu_\pp(\d \theta)\d r \end{equation}    for some strictly positive function  $\psi\in C_b^\infty(\pp M\times [0,\dd])$.
In particular, there exists a constant $c_0>0$ such that
\beq\label{C0} \int_{\pp M\times [0,\dd]} |f(\theta,r)| ^2 \mu_\pp(\d\theta)\d r \le c_0\int_M f^2\d\mu,\ \ f\in \B(M).\end{equation}

\beg{lem}\label{L4}  There exists a constant $c>0$ such that
$$\int_{\pp_r M}\phi^2\d\mu_\pp^r\le c\ss\ll \int_{M}\phi^2\d\mu,\ \ (\ll,\phi)\in \Eig(L):=\Eig_N(L)\cup \Eig_D(L), r\in [0,\dd].$$
 \end{lem}

\beg{proof}  By the symmetry, we only prove the inequality for $r\in [0,\dd/2]$. For $r \in [0,\dd/2]$,   define
$$\gg(s)= \ff{\dd-r}\pi \sin\Big(\ff{(s-r)\pi}{\dd-r}\Big),\ \ s\in [r,\dd].$$ Then $|\nn \gg(\rr_\pp)|\le 1$ and
$$\sup_{M_{r,\dd}} |L\gg(\rr_\pp)| \le \sup_{M_{r,\dd}}|L\rr_\pp| +\ff\pi{\dd-r}\le \sup_{M_{r,\dd}}|L\rr_\pp| +\ff{2\pi}\dd =:c_0<\infty.$$
Let   $N$ be the inward normal unit vector field of $\pp M_{r,\dd}.$ Then
$$N \gg(\rr_\pp)|_{\pp M_{r,\dd}}= \gg'(r)1_{\{\rr_\pp=r\}}- \gg'(\dd)1_{\{\rr_\pp=r\}}=1_{\pp M_{r,\dd}}.$$
So, by integration by parts, there exists a constant $c>0$ such that for any $(\phi,\ll)\in \Eig(L)$,
\beg{align*} &\int_{\pp_r M}\phi^2 \d\mu_\pp^r+  \int_{\pp_\dd M}\phi^2 \d\mu_\pp^\dd=\int_{\pp M_{r,\dd}} \phi^2 N\gg(\rr_\pp) \d\mu_\pp^{r,\dd}\\
&= -\int_{M_{r,\dd}} \big\{\phi^2  L \gg(\rr_\pp) +\<\nn\gg(\rr), \nn\phi^2\>\big\}\d\mu\\
&\le   \int_M(c\phi^2+|\phi|\cdot |\nn\phi|)\d\mu \le c+ \ss{\mu(\phi^2)\mu(|\nn \phi|^2)}= c+\ss\ll.  \end{align*}
\end{proof}
Combining \eqref{4.3} with  Lemmas \ref{L3} and \ref{L4}, we conclude that
\beq\label{HE33}  |\mu(\|\Hess_\phi\|_{HS}^2)-\ll^2|\le c \ll^{\ff 3 2},\ \ (\ll,\phi)\in \Eig_N(L).\end{equation}

\beg{lem}\label{L5}  For any $(\phi,\ll)\in \Eig_N(L)$, under the Fermi coordinates $(\theta,\ll)\in \pp M\times [0,\dd]$ let
\beq\label{H0} h(r)=\ff 1 \ll\int_{\pp M} \phi^2(\theta,r)  \mu_\pp(\d\theta),\ \ r\in [0,\dd].\end{equation} Then there exists a constant $C>0$ depending only on $M_{0,\dd}$ and $L$ such that for $r\in [0,\dd]$ with $h'(r)\ge 0$,
$$h''(r)\ge \ff{|h'(r)|^2}{h(r)}-C.$$
 \end{lem}
\beg{proof} (1) 
 Obviously, we have $h'(r)=\ff 2 \ll \int_{\pp M} \phi\pp_r \phi \d\mu_\pp $ and
\beq\label{HDD}   h''(r) = \ff 2 \ll \int_{\pp M} \big\{(\pp_r\phi)^2 +\phi \pp_r^2\phi\big\}(\cdot,r)\d\mu_\pp, \  \ r\in [0,\dd]. \end{equation}
Let $N_{0,r}$ be the inward unit normal vector field of $\pp M_{0,r}=\{\rr_\pp=0\}\cup\{\rr_\pp =r\}$.
Noting that $\pp_r\phi|_{r=0}= N\phi|_{\pp M}=0$, and $L\phi=-\ll\phi$ implies
$$\phi\psi^{-1} [\pp_r, L] \phi= (L\phi)\psi^{-1} \pp_r \phi -\phi \psi^{-1} L(\pp_r\phi),$$ where $[\pp_r,L]:= \pp_r L- L \pp_r$ is a continuous second order differential operator on the compact set $\{\rr_\pp\le \dd\}$,
by \eqref{**0} and  the integration by parts formula, we have
\beq\label{LPQ0} \beg{split} &\ff 2 \ll \int_{\pp M} \Big(\{\phi \pp_r^2\phi\}(\cdot, r)+ \big\{\phi N^2 \phi\}(\cdot,0)\Big)\d\mu_\pp = \ff 2 \ll  \int_{\pp M_{0,r}}  \big\{\phi \psi^{-1} N_{0,r}\pp_r\phi\big\}\, \d\mu_\pp^{0,r}\\
&= \ff 2 \ll \int_{M_{0,r}} \big(\phi \psi^{-1} L(\pp_r\phi) -\{L(\phi\psi^{-1})\}\pp_r \phi\big)\d\mu +\ff 2 \ll \int_{\pp_r M} (\pp_r \phi) \pp_r(\phi\psi^{-1})\d\mu_\pp^r\\
&\ge \ff 2 \ll \int_{M_{0,r}} \Big(\phi\psi^{-1} [L,\pp_r]\phi - \phi (\pp_r\phi) L\psi^{-1} - 2 \<\nn\phi,\nn\psi^{-1}\> \pp_r \phi\Big)\d\mu \\
&\qquad +\ff 2 \ll \int_{\pp_r M}\psi  (\pp_r \phi) \pp_r(\phi\psi^{-1})\d\mu_\pp \\
&=: I_1+I_2.\end{split}\end{equation}
By \eqref{LPQ} and \eqref{HE33}, we obtain
\beq\label{LPQ1} \ff 2 \ll \bigg|\int_{\pp M}  \{\phi N^2 \phi\}(\cdot,0)\Big)\d\mu_\pp\bigg|\le \aa_1\end{equation}
for some constant $\aa_1>0$. Next, since $[\pp_r, L]$ is a continuous second order differential operator on the compact domain $M_{0,\dd}$, $\mu(\phi^2)=1, \mu(|\nn \phi|^2 ) =\ll$ and \eqref{HE33} holds, we may find out   constants $\aa_2,\aa_3>0$ such that
\beq\label{LPQ2}|I_1|\le \ff {\aa_2}\ll \int_{M} \big\{|\nn \phi|^2+|\phi|(|\nn \phi| + \|\Hess_\phi\|_{HS})\big\} \d\mu\le  \aa_3.\end{equation}
Moreover, obviously
$$I_2= \ff 2 \ll \int_{\pp M}  (\pp_r\phi)^2 \d\mu_\pp+\ff 2 \ll \int_{\pp_r M}    (\phi\pp_r\phi)\pp_r\psi^{-1}\d\mu_\pp^r.$$   Combining this with \eqref{LPQ0}-\eqref{LPQ2}, we we arrive at
=\beq\label{YPP} h''(r)\ge \ff {4}{\ll}\int_{\pp M} (\pp_r\phi)^2\d\mu_\pp  -\aa_1  +\ff 2 \ll \int_{\pp_r M}    (\phi\pp_r\phi)(\pp_r\psi^{-1})\d\mu_\pp^r,\  \ r\in [0,\dd].\end{equation}
Since $h'(r)= \ff 2 \ll \int_{\pp M} \phi\pp_r\phi \d\mu_\pp$, by Cauchy-Schwarz inequality we obtain
$$\ff {4}{\ll}\int_{\pp M} (\pp_r\phi)^2\d\mu_\pp\ge \ff{ |h'(r)|^2}{ h(r)}.$$ This together with \eqref{YPP} implies
\beq\label{YPP2} h''(r)\ge  \ff{ |h'(r)|^2}{ h(r)} -\aa_1  + \ff 2 \ll \int_{\pp_r M}    (\phi\pp_r\phi)(\pp_r\psi^{-1})\d\mu_\pp^r,\  \ r\in [0,\dd].\end{equation}

(2) Since $N\phi^2|_{\pp M}=0$, $ \psi^{-1}\in C_b^\infty(M_{0,\dd})$, and
\beq\label{APP0} \int_{\pp M}(|L\phi^2|+ |\nn \phi^2|)\le 4\ll + 2\ss\ll,\end{equation}  by the integration by parts formula, there exist  constants $\aa_5,\aa_6>0$ such that
\beg{align*} &\ff 2 \ll \int_{\pp_r M}    (\phi\pp_r\phi)(\pp_r\psi^{-1})\d\mu_\pp^r = \ff 1 \ll \int_{\pp_r M} \big\{\pp_r (\phi^2\pp_r \psi^{-1}) - \phi^2 \pp_r^2\psi^{-1}\big\}\d\mu_\pp^r\\
&\ge -\ff 2 \ll \int_{M_{0,r} } |L(\phi^2 \pp_r\psi^{-1})|\d \mu -\aa_5 h(r) \ge -\aa_6  -\aa_5 h(r).\end{align*} This and  \eqref{YPP2} yield
$$ h''(r)\ge \ff{|h'(r)|^2}{h(r)} -\aa_1-\aa_6-\aa_5h(r),\ \ r\in [0,\dd].$$ So, it suffices to find out a constant $c>0$ depending only on $L$ and $M_{0,\dd}$ such that
\beq\label{SPP} h(r)\le c,\ \ r\in [0,\dd].\end{equation}
By \eqref{C0} and $\mu(\phi^2)=1$, we have
$$\int_0^\dd h(r) \d r  = \ff 2 \ll \int_{\pp M\times [0,\dd]}  \phi^2(\theta,r)^2\mu_\pp(\d \theta)\d r \le \ff{2c_0}\ll.$$
So, there exists $r_0\in [0,\dd]$ such that
\beq\label{SPP1} h(r_0)\le \ff{2c_0}{\dd\ll}\le \ff{2c_0}{\dd\ll_1}:=c_1,\end{equation} where $\ll_1$ is the first non-trivial Neumann eigenvalue of $L$ on $M$. On the other hand, by Integration by parts formula and noting that $\pp_r\phi(\theta,r)=0$ for $r=0$, for $r\in [0,\dd]$ we have 
\beg{align*} |h'(r)| &= \ff 1 \ll \bigg|\int_{\pp_r M} (\pp_r\phi^2)\psi^{-1}\d\mu_\pp^r\bigg| = \ff 1 \ll \bigg|\int_{\pp M_{0,r}} (\pp_r\phi^2)\psi^{-1}\d\mu_\pp^{0,r}\bigg|\\
&\le \ff 1 \ll  \int_{M_{0,r}} \big(|\psi^{-1} L\phi^2| +|\<\nn\psi^{-1},\nn \phi^2\>|\big)\d\mu.\end{align*}
Combining this with \eqref{APP0}, we find out a constant $c_2>0$ such that 
$$|h'(r)| \le c_2,\ \ r\in [0,\dd].$$
This together with \eqref{SPP1} implies \eqref{SPP} for $c=    c_1 +\dd c_2.$ Then the proof is finished.
 \end{proof}

\beg{proof}[Proof of Lemma \ref{TP}] Due to  Lemma \ref{L5}, this result can be proved by modifying the argument in    \cite[Proof of Proposition 2.4]{BHT}.
Let $C$ be the constant in Lemma \ref{L5}.

(1) We first prove that for large enough $\ll>0$,
\beq\label{*D4} r\in [0,2\dd/3] \ \text{with}\   h'(r)>0 \ \text{implies}\  |h'(r)|^2 < 4C h(r).\end{equation}
If the assertion is not true, then there exists $r_0\in [0, 2\dd/3]$ such that $h'(r_0)>0$ and $|h'(r_0)|^2\ge 4C h(r_0)(1+h(r_0))$. Then
  by Lemma \ref{L5}
\beq\label{*D'}\beg{split}  &\ff{\d}{\d r} \{|h'(r)|^2-4 C h(r)\}\big|_{r=r_0}= 2 h'(r_0)h''(r_0) - 4C h'(r_0) \\
&\ge \ff{2(h'(r_0))^3}{h(r_0)}- 6Ch'(r_0) = 2 C h'(r_0)\ge  4 C\ss{Ch(r_0)}>0.\end{split} \end{equation} So, there exists $\vv\in (0, \dd-r_0)$ such that   $h'(r)>0$ and $|h'(r)|^2\ge 4 C h(r)$ hold for $r\in [r_0,r_0+\vv]$. By a continuity argument we  conclude that
\beq\label{*D3} h'(r)>0,\   |h'(r)|^2\ge 4C h(r),\ \ r\in [r_0,\dd].\end{equation}
Indeed, if not then
$$r_1:=\inf\{r\in [r_0, \dd]:     |h'(r)|^2< 4 C h(r)\}\in [r_0+\vv,\dd]\subset (r_0,\dd].$$  We have $h'(r)>0$  for $r\in [r_0,r_1]$ and $|h'(r_1)|^2= 4 C h(r_1)$, so that \eqref{*D'}  holds for $r_1$ replacing $r_0$. Thus, due to continuity, there exists $r_2\in [r_0, r_1)$ such that
$$\ff{\d}{\d r} \{|h'(r)|^2- 4Ch(r)\} >0,\ \ r\in [r_2, r_1].$$
Since by the definition of $r_1$ we have $|h'(r_2)|^2- 4C h(r_2)\ge 0,$   this implies
$$|h'(r_1)|^2 -4Ch(r_1)=\sup_{r\in [r_2,r_1]} \{|h'(r)|^2 -4C h(r)\} >0,$$ which contradicts to $|h'(r_1)|^2 = 4C h(r_1).$  So, \eqref{*D3} holds and thus,
$$\ff{\d}{\d r} \ss {h(r)}\ge \ss C,
 \ \ r\in [r_0,\dd].$$ This implies $h(r)\ge C(r-r_0)^2$ for $r\in [r_0,\dd]$, and hence, by \eqref{C0},
 \beg{align*} &\ff{C} 3  \Big(\dd- \ff{2\dd}3\Big)^3\le  C\int_{r_0}^\dd (r-r_0)^2\d r \\
  &\le  \int_0^\dd h(r)\d r \le \ff{c_0}\ll \int_M\phi^2\d\mu =\ff{c_0}\ll,\end{align*} which is impossible for large enough $\ll$. The contradiction means that for large enough $\ll>0$, \eqref{*D4} holds.

 (2) We then prove that for large $\ll>0$,
 \beq\label{*D5} |h'(r)|^2\le 5C h(r),\ \ r\in [0,2\dd/3].\end{equation}
 By the Neumann condition we have $h'(0)=0$, so that the inequality in \eqref{*D5} holds in a neighborhood of $0$. Thus, if \eqref{*D5} does not hold, then
 $$r_2:= \inf\big\{r\in [0,2\dd/3]: |h'(r)|^2 >5C h(r)\big\}\in (0, 2\dd/3]$$ exists, and
 $$|h'(r_2)|^2 = 5C h(r_2),\ \ \ff{\d}{\d r} \big\{|h'(r)|^2 - 5C h(r)\big\}\big|_{r=r_2} \ge 0.$$
 Combining this with \eqref{*D4}, we obtain $h'(r_2)< 0$ and
 $$0 \le \ff{\d}{\d r} \big\{|h'(r)|^2 - 5C h(r)\big\}\big|_{r=r_2} = 2 h'(r_2)h''(r_2) -5Ch'(r_2) = 2 h'(r_2)\Big(h''(r_2)-\ff{5 C}2\Big).$$
 So, $h''(r_2)\le \ff{5C} 2.$ But by Lemma \ref{L5} and $|h'(r_2)|^2=5Ch(r_2)$ we have
 $$h''(r_2)\ge \ff{|h'(r_2)|^2}{h(r_2)} -C= 4C,$$ which is a contradiction.  Therefore, \eqref{*D5} has to be true.

 (3) By \eqref{*D5}, when $\ll>0$ is large enough  we have
 $$\ff{\d }{\d r} \ss{h(r)} = \ff{h'(r)} {2\ss{h(r)} } \ge -\ff{\ss{5Ch(r)}}{2\ss{h(r)}}=-\ff{\ss{5C}}2 =:-c,\ \ r\in [0, 2\dd/3].$$
 So,
 $\ss{h(r)}\ge \ss{h(0)} -\ss c\, r$ holds for $r\in [0,2\dd/3].$ Let $M= \ff{16 (1\lor c_0)}{c}$, where $c_0$ is in \eqref{C0}. If $h(0)\ge Mc\ll^{-2/3}$, we would have
 $$\ss{h(r)}\ge \ss{Mc} (\ll^{-1/3}-r) \ge \ff{\ll^{-1/3}}2\ss{Mc},\ \ r\in [0, \ll^{-1/3}/2],$$ where we take $\ll>0$ large enough such that $\ll^{-1/3}/2\le 2\dd/3.$
 Combining this with \eqref{C0} and \eqref{H0}, we arrive at
 \beg{align*} \ff{2c_0}\ll &\le \ff{Mc}{8\ll} \le \int_0^{\ll^{-1/3}/2} h(r)\d r \le \ff 1 \ll \int_{\pp M\times [0,\dd]} \phi(\theta,r)^2 \mu_\pp(\d\theta)\d r\\
 &\le \ff{c_0}\ll \int_M\phi^2\d\mu =\ff{c_0}\ll\end{align*} for large enough $\ll>0$, which is however impossible. This means that when $\ll>0$ is large enough we have
 $h(0)\le Mc \ll^{-2/3},$ equivalently,
 $$\int_{\pp M}\phi^2\d\mu_\pp \le Mc \ll^{1/3},
 $$ which completes the proof.
\end{proof}

We are now ready to prove Theorem \ref{T4.2}.

\beg{proof}[Proof of Theorem \ref{T4.2}] Since $L\phi=-\ll\phi$, estimate \eqref{4.3} follows from \eqref{IF1}. Moreover, estimates in \eqref{4C20} are included in Lemma \ref{L3} and Lemma \ref{TP}. Combining \eqref{4.3} with \eqref{4C20} we prove \eqref{4.2}.
\end{proof}






\end{document}